\documentclass[11pt,twoside]{amsart}
\usepackage[margin=1.2in]{geometry}
\usepackage[utf8]{inputenc}
\usepackage{amsmath}
\usepackage{quiver}
\usepackage{mathtools}
\usepackage{graphicx} 
\graphicspath{ {./images/} }
\usepackage{wrapfig}
\usepackage{amsmath}
\usepackage{amsthm}
\usepackage{amsfonts}
\usepackage{amssymb}
\usepackage[all]{xy}
\usepackage{alltt}
\usepackage{enumitem}
\usepackage{stmaryrd}
\usepackage{comment}
\usepackage{xspace}
\usepackage{xcolor}
\usepackage{comment}
\usepackage{float}

\theoremstyle{plain}
\newtheorem{prop}{Proposition}[section]

\newtheorem{lem}[prop]{Lemma}

\newtheorem{corollary}[prop]{Corollary}

\newtheorem{thm}[prop]{Theorem}
\theoremstyle{definition}

\newtheorem{cor}[prop]{Corollary}
\newtheorem{remark}[prop]{Remark}
\newtheorem{example}[prop]{Example}
\newcommand{\expl}{\noindent{\bf{Example}}}
\newcommand{\rmrk}{\noindent{\bf{Remark }}}
\newcommand{\zed}{\mathbb{Z}}
\newcommand{\field}{\mathbb{F}}
\renewcommand{\qedsymbol}{$\blacksquare$}
\newcommand{\trank}{\ensuremath{\mathrm{d}}}
\newcommand{\pf}{\textbf{Proof}.\quad}
\renewcommand{\brack}[1]{\langle #1 \rangle}

\usepackage[
pdfauthor={ESYZ},
pdftitle={},
pdfstartview=XYZ,
bookmarks=true,
colorlinks=true,
linkcolor=blue,
urlcolor=blue,
citecolor=blue,
bookmarks=false,
linktocpage=true,
hyperindex=true
]{hyperref}

\title{ Rings with an elementary abelian $p$-group of units}
\author{Sunil Chebolu}

\address{Department of Mathematics \\
Illinois State University \\
Normal, IL 61790, USA}
\email{schebol@ilstu.edu, jacorry@ilstu.edu, evgrimm@ilstu.edu, abhatfi@ilstu.edu}

\author{Jeremy Corry}
\author{Elizabeth Grimm}
 \author{Andrew  Hatfield}

\date{\today}

\begin{document}

\thanks{The first author is supported by Simons Foundation: Collaboration Grant for Mathematicians (516354). }


\begin{abstract}
 What are all rings $R$ for which $R^\times$ (the group of invertible elements of $R$ under multiplication) is an elementary abelian $p$-group? We answer this question for finite-dimensional commutative $k$-algebras, finite commutative rings,  modular group algebras, and path algebras. Two interesting byproducts of this work are a characterization of Mersenne primes and a connection to Dedekind's problem.
\end{abstract}
\maketitle


\section{Introduction} 

 Robert Gilmer \cite{GIL} classified all finite commutative rings $R$ with identity such that $R^{\times}$, the group of units of $R$, is cyclic. He showed that these rings are isomorphic to a finite product of rings $R_i$ from the list below such that for $i \ne j$, $\gcd(|R_i^\times|, |R_j^\times|) = 1$:

\begin{itemize} 
    \item The finite field $\mathbb{F}_{p^k}$,
    \item $\mathbb{Z}/p^m$ where $p$ is an odd prime and $m > 1$,
    \item $\mathbb{Z}/4$,
    \item $\mathbb{F}_p[x]/(x^2)$ where $p$ is any prime,
    \item  $\mathbb{F}_2[x]/(x^3)$, and
    \item $\mathbb{Z}[x]/(4, 2x, x^2-2)$.
\end{itemize}

The unit groups of the above rings are well-known. Therefore, given any positive integer $m$, using Gilmer's theorem, one can write down all finite commutative rings whose unit group is $C_m$, the cyclic group of order $m$.

 After cyclic groups, an important family of finite groups is the class of elementary abelian $p$-groups - a direct sum of copies of the cyclic group $C_p$, where $p$ is a fixed prime. These groups feature in many areas, including topology, group cohomology, and representation theory. This paper parallels Gilmer's result by classifying all finite commutative rings whose unit group is an elementary abelian $p$-group. In addition, we also obtain a classification of all finite-dimensional $k$-algebras, modular group algebras, and path algebras for which the unit group is an elementary abelian $p$-group.

The first author  studied this problem for the families $\mathbb{Z}_n$ \cite{24}, $\mathbb{Z}_n[x_1, x_2, \ldots, x_m]$ \cite{12}, and group algebras over fields \cite{CLY} under the name of diagonal property. A ring $R$ is said to have the diagonal property if 1's in the multiplication table for $R$ fall only on the main diagonal. For any prime $p$, a ring $R$ is  said to be a $\Delta_p$-ring if all units $u$ in $R$ satisfy $u^p = 1$. It is easy to see that $R$ has the diagonal property if and only if it is a $\Delta_2$-ring. Furthermore, these conditions are equivalent to $R^\times$ being an elementary abelian $2$-group; see Proposition \ref{elemabelian}.  We now state our main results.

\begin{thm}Classifying rings whose unit group is an elementary abelian $2$-group.
\begin{enumerate} 
\item Let $R$ be a finite-dimensional commutative $k$-algebra where $k$ is a field. $R$ is a $\Delta_2$-ring if and only if $R$ is isomorphic to a finite product of rings that are either $\mathbb{F}_3$ or a quotient of a truncated polynomial ring of the form $\mathbb{F}_2[x_1, x_2, \ldots, x_l]/(x_1^2, x_2^2, \ldots, x_l^2)$, $l \ge 1$.
\item Let $R$ be a finite commutative ring. $R$ is a $\Delta_2$-ring if and only if $R$ is isomorphic to any quotient of a finite product of rings $R_i$, where each $R_i$ is any of the following rings:
\begin{itemize}
    \item $\mathbb{F}_2[x_1, x_2, \ldots, x_l]/(x_1^2, x_2^2, \ldots, x_l^2)$, $l \ge 1$.
    \item $\mathbb{F}_3$
    \item $\mathbb{Z}_4[x_1, x_2, \ldots, x_l]/( r(r+2)  \colon r \in \mathcal{P} )$ where $\mathcal{P} = (2, x_1 -1, x_2-1,\ldots, x_l-1) \subseteq  \mathbb{Z}_4[x_1, x_2, \ldots, x_l]$.
    \item  $\mathbb{Z}_8[x_1, x_2, \ldots, x_l]/( r(r+2)  \colon r \in \mathcal{P} )$ where $\mathcal{P} = (2, x_1 -1, x_2-1,\ldots, x_l-1) \subseteq \mathbb{Z}_8[x_1, x_2, \ldots, x_l]$.
\end{itemize}
\item Let $n$ be a positive integer and $G$ be a finite group. $\mathbb{Z}_nG$ is a $\Delta_2$-ring if and only if $n = 2, 3, 6$ and $G$ is an elementary abelian $2$-group or $n = 4, 12$ and $G = C_2$.
\item Let $Q$ be a finite and acyclic quiver. The path algebra $kQ$ is a $\Delta_2$-ring if and only if $k=\mathbb{F}_3$ and $Q$ is trivial   or $k=\mathbb{F}_2$ and $Q$ has no directed path of length $2$.
\end{enumerate}
\end{thm}

Analyzing the quotients mentioned in the above classification led us to Dedekind's problem, which is wide open; see Remark \ref{Dedekind}. Also, note that the so-called exceptional primes (2 and 3) feature in parts 1 and 4 of the above theorem. These primes play a special role in many places. For instance, they appear in work on generating hypothesis for the stable module categories. Combining this theorem with earlier work done in Tate cohomology \cite[Theorem 1.1]{GH} gives the following intriguing characterization of these exceptional primes.

\begin{cor}
Let $p$ be a prime number. Then the following are equivalent.
\begin{enumerate}
    \item $p = 2$ or $3$
    \item There is a finite-dimensional $\mathbb{F}_p$-algebra whose unit group is an elementary abelian $2$-group.
    \item The Tate cohomology functor $\hat{H}(C_p, -)$ is faithful on the stable module category of finitely generated $\mathbb{F}_pC_p$-modules.
\end{enumerate}
\end{cor}

When looking for rings whose unit group is an elementary abelian $p$-group, with $p$ odd, it is natural to ask which elementary abelian $p$-groups occur as unit groups of rings. That is Fuchs' problem for elementary abelian $p$-groups and was answered in \cite{CL}. An abelian $p$-group occurs as a unit group of a ring if and only if $p=2$ or a Mersenne prime (a prime of the form $2^m-1$ for some $m$). That is the reason why Mersenne primes feature in the following theorem.

Moreover, for odd primes $p$, the $\Delta_p$ condition  and $R^\times$ being an elementary abelian $p$-group are equivalent for commutative rings;  see Proposition \ref{elemabelian}.

\begin{thm}Classifying rings whose unit group is an elementary abelian $p$-group, where $p$ is an odd prime.
\begin{enumerate} 
\item Let $R$ be a finite-dimensional commutative $k$-algebra. $R$ is a $\Delta_p$-ring if and only if  $p$ is a Mersenne prime and  $R = (\field_2)^a \times (\mathbb{F}_{p+1})^b$ for some nonnegative integers $a$ and $b$.
\item A finite commutative ring $R$ is is a $\Delta_p$-ring if and only if $p$ is  Mersenne prime and  $R$ is isomorphic to $R = (\field_2)^a \times (\mathbb{F}_{p+1})^b$ for some nonnegative integers $a$ and $b$.
\item  $\mathbb{Z}_nG$ ($G$ abelian) is a $\Delta_p$-ring if and only if $n=2$, $G$ is an elementary abelian $p$-group, and $p$ is some Mersenne prime.
\item Let $Q$ be a finite acyclic quiver. The path algebra $kQ$ is $\Delta_p$ if and only if $k = \mathbb{F}_2$ or $\mathbb{F}_{p+1}$ where $p$ is a Mersenne prime and $Q$ is trivial.

\end{enumerate}
\end{thm}
 It is worth noting that the above theorem gives some characterizations of Mersenne primes; see Theorem 1.1 \cite{CLY} for similar characterizations. 

\begin{cor} Let $p$ be an odd prime. Then the following are equivalent. 
\begin{enumerate}
\item $p$ is a Mersenne prime.
\item There exists a  finite-dimensional commutative  $k$-algebra that is $\Delta_p$.
\item There exists a finite commutative  $\Delta_p$-ring. 
\item There exists a finite acyclic quiver $Q$ and a field $k$ such that the path algebra $kQ$ is $\Delta_p$-ring.
\end{enumerate}

\end{cor}

We see from the above results that $\Delta_p$-rings are rare. For instance,  working over a field $k$, it is easy to see that for $n > 1$,  $M_n(k)$ is never a $\Delta_p$-ring.

The paper is organized as follows. We begin in Section \ref{preliminaries} with some preliminaries. We then classify the $\Delta_p$-rings for finite-dimensional commutative $k$-algebras (Section \ref{fdcka}), finite commutative rings (Section \ref{fcr}), modular group algebras (Section \ref{mgr}), and path algebras of quivers (Section \ref{pa}).

\vskip 5mm
\noindent
\textbf{Acknowledgements:}
We thank Dave Benson, Jon Carlson, Srikanth Iyengar, and  Richard Stanley for the discussions related to this paper. This led to a better understanding of the lattice of ideals in truncated polynomial rings and its relation to the Dedekind problem; see Remark \ref{Dedekind}. We also thank an anonymous referee for several comments that improved the exposition.

\section{Preliminaries} \label{preliminaries}

This section collects some background material and lemmas we need to prove our main results.   We begin with the relationship between $\Delta_p$ and $R^\times$ being an elementary abelian $p$-group.

\begin{prop}\label{elemabelian} Let $R$ be a unital ring.  Consider the following statements.
\begin{enumerate}
 \item $R$ is a $\Delta_p$-ring.
    \item  $R$ is a ring such that $R^\times$ is an elementary abelian $p$-group
   
\end{enumerate}
When $p = 2$ the above two statements are equivalent. When $p > 2$, $2$ implies $1$, and $1$ implies $2$ provided $R^\times$ is commutative. 
\end{prop}

\begin{proof}
 $(2) \implies (1)$ with $p=2$ is the only non-trivial part, and that follows from an exercise in group theory: any group $G$ such that $g^2 =e$ for all $g$ in $G$ is abelian.
If $G$ is  abelian such that $g^2 =e$ for all $g$ in $G$, then $G$ is a vector space over $\field_2$. If we let
 $r$ be the dimension of $G/\field_2$, then we see that $G \cong \field_2^r \cong C_2^r$,
 so $G$ is an elementary abelian 2-group.
\end{proof}

We now summarize results from \cite{24, 12, CLY} related to the  $\Delta_p$ condition.

\begin{thm} Examples of $\Delta_2$ and $\Delta_p$ rings.
\begin{enumerate}
    \item \cite{24} $\mathbb{Z}_n$ is a $\Delta_2$ ring if and only if $n$ divides 24.
    \item \cite{12} $\mathbb{Z}_n[x_1, \ldots, x_m]$ is a $\Delta_2$ ring if and only if $n$ divides 12 and $m \ge 1$.
    \item \cite{CLY} For a field $k$ and group $G$, $kG$ is $\Delta_2$ if and only if $k =\field_2$ or $\field_3$ and $G$ is an elementary abelian $2$-group of possibly infinite rank.
    \item \cite{CLY} Let $p$ be an odd prime and let $k$ be a field, and $G$ be an abelian group. $kG$ is $\Delta_p$ if and only $p$ is Mersenne and $kG$ is either $\field_2C_p^r$ or $\field_{p+1}C_p^r$ where $0 < r \le \infty$.
\end{enumerate}
\end{thm}

\begin{lem} \label{basicproperty}
A subring of a $\Delta_p$-ring is again a $\Delta_p$-ring. A direct product of rings is $\Delta_p$  if and only if each factor is a $\Delta_p$-ring.
\end{lem}

\begin{lem} \label{24}
Let $R$ be a ring of characteristic $n\, (>1)$. If $R$ is a $\Delta_2$-ring, then $n$ divides $24$.
\end{lem}
\begin{proof}
 For a $\Delta_2$-ring $R$ with characteristic $n$, $\mathbb{Z}_n$ is a subring of $R$, and therefore a $\Delta_2$-ring. It was shown in \cite{24} that $\mathbb{Z}_n$ is a $\Delta_2$-ring if and only if $n$ divides $24$.
\end{proof}

\begin{lem}\label{deltapfield} \cite{CLY}
    Let $p$ be a prime. A field $k$ is $\Delta_p$ if and only if $k = \field_2$ or $k = \field_3$ with 
    $p = 2$, or $k = \field_2 $ or $\field_{p+1}$ with $p$ a Mersenne prime. In particular, these conditions on $k$ hold for any $k$-algebra that is $\Delta_p$.
\end{lem}

An Artinian ring is a ring that satisfies the descending chain condition on its ideals. That is, for every descending chain of ideals $I_1 \supseteq I_2 \supseteq  \ldots, \supseteq  I_k \supseteq \ldots, $ in $R$, there is an integer $k$ such that  $I_n = I_{n+1}$ for all $n > k$. Recall that finite commutative rings and finite-dimensional $k$-algebras are Artinian. Artinian rings have the following Artin decomposition \cite{Atiyah}: every Artinian ring $R$ is a direct product of Artin local rings $R_i$,
\[ R = R_1 \times R_2 \times \ldots \times R_n.\]
Moreover,  an Artinian local ring has a unique prime ideal. The following lemma will be used in our analysis.

\begin{lem} \label{keylemma}
Let $R$ be an Artinian ring with a unique prime ideal $P$. Then $P$ is the set of all nilpotent elements, and $R \setminus P$ is the set of all units. Moreover, $R$ is generated as a ring by its units. In particular, when $R^\times$ is an elementary abelian $p$-group of rank $t$, $R$ is a quotient of the group ring $\mathbb{Z}_n[C_p^t]$, where $n$ is the characteristic of $R$.
\end{lem}

\begin{proof}
The first statement is well-known. To see that $R \setminus P$ generates $R$, note that for every nilpotent $\eta$, the element $u:= 1+\eta$ is a unit. This shows that $\eta = u - 1$, proving that the units generate all elements of the ring. Let $n \,(\ge 0)$ be the characteristic of $R$. When $R^\times = R \setminus P \cong C_p^t$ generates $R$, there is a surjective ring homomorphism 
\[ \mathbb{Z}_n[C_p^t] \rightarrow R.\]
This shows that $R$ is a quotient of the group ring $\mathbb{Z}_n[C_p^t]$.
\end{proof}

\begin{thm}[Artin-Wedderburn]\label{WA}
    Any semisimple ring $R$ is isomorphic to a product of finitely many matrix 
    rings over division rings $D_i$ for some dimension $n_i$, both of which are uniquely determined
    by permutation of $i$, i.e. 
    \begin{align*}
        R = \prod M_{n_i}(D_i).
    \end{align*}
\end{thm}

An important family of semisimple rings comes from Maschke's theorem, which states that the group algebra $kG$ is semisimple when the characteristic of the field $k$ is relatively prime to the order of the finite group $G$. 

\begin{cor} \label{commutativeSS}
Every finite commutative semisimple ring is a direct product of finite fields. In particular, any quotient $\field_p G /\!\!\sim$ of a group algebra $kG$, where $G$ is a group whose order is relatively prime to $p$, is a finite product of finite fields of characteristic $p$.
\end{cor}

\begin{proof}
By the Artin-Wedderburn theorem, every semisimple ring $R$ is a direct product of matrix rings over division rings. If $R$ is commutative, these matrix rings  have order $1 \times 1$. That is, $R$ is a product of division rings. Furthermore, by Wedderburn's little theorem, we know that finite division rings are fields. This shows that $R$ is a product of finite fields. When $p$ and $G$ are as given, Maschke's theorem implies that $\field_p G$ is a product of finite fields of characteristic $p$. Since a quotient of a product of fields is again a product (with possibly fewer factors), the second statement follows.
\end{proof}

\section{Finite dimensional commutative $k$-algebras} \label{fdcka}

We want to classify all finite-dimensional $k$-algebras that are $\Delta_p$.  We begin with case $p=2$. 
For any ring, $A$, $A/\!\!\sim$  will denote an arbitrary quotient of  $A$, and
$\prod R_l$ denotes a finite product of rings of the form $R_l$.

\begin{thm}
Let $R$ be a finite-dimensional commutative $k$-algebra. Then $R$ is $\Delta_2$ if and only if $R \cong \prod \mathbb{F}_3$ or a quotient of $\prod \mathbb{F}_2[x_1,...,x_l]/(x_1^2,...,x_l^2)$.
\end{thm}

\begin{proof}
Let $R$ be as given and assume that it is $\Delta_2$. Then $R$ has an Artin decomposition:
\[R = R_1 \times \ldots \times R_l,\]
where each $R_i$ is a local ring. By Lemma $\ref{basicproperty}$, we know that $R$ is $\Delta_2$ if and only if each $R_i$ is $\Delta_2$. So, we assume that $R$ is a finite-dimensional local  $k$-algebra. Moreover, by Lemma \ref{deltapfield}, it is enough to assume that $k = \field_2$ or $\field_3$. Artin local rings have a unique prime ideal. Therefore, Lemma \ref{keylemma} implies that  $R \cong \mathbb{F}_2[C_2^r]/\!\!\sim$ or $R \cong \mathbb{F}_3[C_2^r]/\!\!\sim$. 

 Let us first consider the case $R \cong \mathbb{F}_2[C_2^r]/\!\!\sim$. Note that 
 \[\mathbb{F}_2[C_2^r] \cong \field_2[x_1, \ldots, x_r] / (x_1^2 - 1, \ldots, x_r^2 - 1) \cong \field_2[t_1, \ldots, t_r] / (t_1^2, \ldots, t_r^2),\]
 The first isomorphism sends a generator in the $i$th factor in $C_2^r$ to $x_i$.   The second isomorphism is obtained using a change of variables: $t_i = x_i-1$.  Note that  $t_i^2 = (x_i -1)^2 = x_i^2 - 2x_i + 1 = x_i^2 - 1$ in a field with characteristic 2. Because $R$ is isomorphic to a quotient of ${F}_2[C_2^r]$, we see that 
 $R \cong \field_2[t_1, \ldots, t_r]/I$, where $ \brack{t_1^2, \ldots, t_r^2 } \subseteq I$. Now let $x$ be a unit in $R$. Using Lemma \ref{keylemma}, we can write $x$ as $\overline{1 + \eta}$ for some  $\eta$ in $(t_1, t_2, \ldots t_r)$.  Then $x^2 = \overline{1^2 +\eta^2} = \overline{1}$ because characteristic of $R$ is $2$ and $\eta^2 = 0$ in $R$. This shows that $R$ is a $\Delta_2$-ring.
 
 The second case is $R \cong \mathbb{F}_3[C_2^r]/\!\!\sim$. In this case, Corollary \ref{commutativeSS} implies that $R \cong \prod  \field_{3^{r_i}}$. Taking units, we see that 
 \[ R^\times \cong \prod  (\field_{3^{r_i}})^\times \cong \prod C_{3^{r_i}-1}.\]
This is an elementary abelian $2$-group if and only if $3^{r_i}-1 = 2$, or $r_i=1$ for all $i$. This shows that $R \cong \prod \field_3$, which is clearly $\Delta_2$.
\end{proof}

\begin{remark}  \label{Dedekind} Note that in the above theorem, we classified the finite-dimensional $\field_2$-algebras that are $\Delta_2$ in terms of the quotients of the ring $\mathbb{F}_2[x_1,...,x_l]/(x_1^2,...,x_l^2)$. So, a natural question is: what are all quotients of this ring?

The quotients of this ring correspond to the ideals of the polynomial ring  $\mathbb{F}_2[x_1,...,x_l]$ that contain the ideal $(x_1^2,...,x_l^2)$. These ideals form a partially ordered set under inclusion. 
A complete description of this lattice seems hopeless. Even the sublattice of ideals generated by sets of monomials in $\mathbb{F}_2[x_1,...,x_l]$ is not well-understood. The latter is the free distributive lattice (\cite{GG}) $\text{FD}(l)$ on $l$ generators, and computing its cardinality is the Dedekind problem, which is open. The exact values are known only for $1 \le  l \le 8$:  3, 6, 20, 168, 7581, 7828354, 2414682040998, 56130437228687557907788 (sequence A000372 in the OEIS).

Here are the Hasse diagrams for the poset of ideals for $l=2$. The left diagram corresponds to the lattice of all ideals of the quotient ring  $\field_2[x_1, x_2]/(x_1^2, x_2^2)$, the right one to those ideals that are generated by monomials.

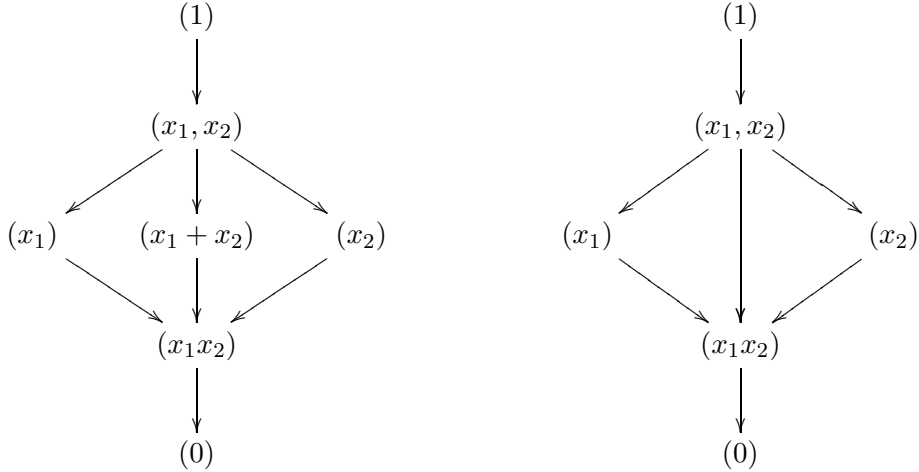
\begin{center}

\begin{minipage}[c]{0.44\textwidth}
\xymatrix{
& (1) \ar[d] & \\
&(x_1, x_2) \ar[dl] \ar[d] \ar[dr] & \\
(x_1) \ar[dr] & (x_1+x_2) \ar[d] & (x_2) \ar[dl]\\
&(x_1x_2) \ar[d] & \\
&(0) &
}
\end{minipage}
\hspace{0.03\textwidth}
\hspace{0.08\textwidth}
\begin{minipage}[c]{0.44\textwidth}
\xymatrix{
& (1) \ar[d] & \\
&(x_1, x_2) \ar[dl] \ar[dd] \ar[dr] & \\
(x_1) \ar[dr] & \ar[d] & (x_2) \ar[dl]\\
&(x_1x_2) \ar[d] & \\
&(0) &
}
\end{minipage}
\end{center}
\end{remark}

The actual number of quotients of the ring $\mathbb{F}_2[x_1,...,x_l]/(x_1^2,...,x_l^2)$ is at least as big as $|\text{FD}(l)|$ and it can be computed using GAP for small values of $l$. This gave $3,7,47,4979$ for $1 \le l \le 4$.

\begin{thm}
 Let $R$ be a finite-dimensional commutative $k$-algebra, and let $p$ be an odd prime. $R$ is a $\Delta_p$-ring if and only if  $p$ is a Mersenne prime and  $R = (\field_2)^a \times (\mathbb{F}_{p+1})^b$ for some nonnegative integers $a$ and $b$.
\end{thm}

\begin{proof}
Let $R$ be as stated in the theorem and assume that $R$ is $\Delta_p$ for some odd prime $p$. As argued in the proof of the above theorem, we may assume that $R$ is a local Artin $k$-algebra. 
Then by Lemma \ref{deltapfield}, $k$ is either $\field_2$ or $\field_{p+1}$ where $p$ is a Mersenne prime.  Then Lemma \ref{keylemma} implies that  that 
 $R \cong \mathbb{F}_2[C_p^r]/\!\!\sim$ or $R \cong \mathbb{F}_{p+1}[C_p^r]/\!\!\sim$.
Taking products of local algebras and using Lemma \ref{commutativeSS}  we get the classification stated in the theorem because the only finite fields that are $\Delta_p$ are $\field_2$ and $\field_{p+1}$ with $p$ Mersenne.
\end{proof}

\begin{cor}
An odd prime $p$ is Mersenne if and only if there exists a finite-dimensional commutative $k$-algebra that is a $\Delta_p$-ring.
\end{cor}

\section{Finite commutative rings} \label{fcr}

Let $R$ be a finite commutative ring. Note that $R$ is Artinian, and therefore it can be written as a product of Artin local rings:  $R = R_1 \times R_2 \ldots \times  R_l$. Recall that $R$ is $\Delta_2$ if and only if $R_i$ is $\Delta_2$ for all $i$. 
So, we may assume that $R$ is finite and has only one prime ideal. 
Let $n$ be the characteristic of $R$. Then Lemma \ref{24} implies that $n$ divides $24$.

If $n=2$ or $3$, then $R$ is a finite-dimensional $\mathbb{F}_2$ or $\mathbb{F}_3$-algebra. By the previous section, we know that $R$ is either a product of $\mathbb{F}_3$'s or a quotient of a  product of rings of the form $\mathbb{F}_2[x_1, x_2, \ldots, x_l]$. By the Chinese Remainder Theorem, every ring of characteristic 6 is a product of two rings, one of characteristic 2 and one of characteristic 3. So this also takes care of rings of characteristic 6.

It is enough to consider rings of characteristics 4 and 8. The $\Delta_2$ rings in these two cases, when combined with products of $\mathbb{F}_3$, will complete the cases of characteristics 12 and 24 using the Chinese Remainder Theorem.

By Lemma \ref{keylemma}, any local Artin ring of characteristic $4$ that is $\Delta_2$ must  be a quotient of 
\[\mathbb{Z}_4[C_2^l] \cong \mathbb{Z}_4[x_1, x_2, \ldots, x_l]/(x_1^2-1, x_2^2-1, \ldots, x_l^2-1).\]

So, what are all the $\Delta_2$-quotients? We begin by determining the unique prime ideal of this ring.

\begin{lem} \label{Rislocal}
The ring $R = \mathbb{Z}_4[x_1, x_2, \ldots, x_l]/(x_1^2-1, x_2^2-1, \ldots, x_l^2-1)$ has only one prime ideal, and that corresponds to $\mathcal{P} := (2, x_1-1, x_2-1, \ldots, x_l-1)$ in $\mathbb{Z}_4[x_1, x_2, \ldots, x_l]$.
\end{lem}

\begin{proof}
 Let $\mathcal{I}$ be any prime ideal of $\mathbb{Z}_4[x_1, x_2, \ldots, x_l]$ that contains $(x_1^2-1, x_2^2-1, \ldots, x_l^2-1)$. Since we are working in characteristic $4$, we have $2^2 = 4 = 0 \in \mathcal{I}$. Since $\mathcal{I}$ is a prime ideal, $2$ belongs to $\mathcal{I}$. For all $i$, we have $x_i^2-1$ and $2$ belong to $\mathcal{I}$. This means $(x_i^2-1)+2= x_i^2+1$, and hence $x_i^2+1-2x_i =(x_i-1)^2$ belong to $\mathcal{I}$. This shows that $x_i-1$ is in $\mathcal{I}$ for all $i$. The ideal generated by $(2, x_1-1, \ldots, x_l-1)$ is maximal in $\mathbb{Z}_4[x_1, x_2, \ldots, x_l]$ because the quotient is the residue field $\field_2$. This shows that $\mathcal{I} = \mathcal{P} =(2, x_1-1, \ldots, x_l-1)$.
\end{proof}

The next proposition characterizes the quotients of $\mathbb{Z}_4[x_1, x_2, \ldots, x_l]/(x_1^2-1, x_2^2-1, \ldots, x_l^2-1)$ that are $\Delta_2$.

\begin{prop}
Let $R = \mathbb{Z}_4[x_1, x_2, \ldots, x_l]/(x_1^2-1, x_2^2-1, \ldots, x_l^2-1)$ and let $\mathcal{P} = (2, x_1-1, x_2-1, \ldots, x_l-1)$. A  quotient ring $S$, of $R$, is a $\Delta_2$-ring if and only if 
$S = \mathbb{Z}_4[x_1, x_2, \ldots, x_l]/\mathcal{J}$ where $\mathcal{J}$ contains $( \eta(\eta +2) \colon \eta \in \mathcal{P} )$.
\end{prop}

\begin{proof}
We begin by noting that any quotient ring of $R$ is of the form 
\[ S = \mathbb{Z}_4[x_1, x_2, \ldots, x_l]/\mathcal{J},\] 
where $(x_1^2-1, \ldots, x_l^2-1) \subseteq  \mathcal{J} \subseteq \mathcal{P}$. (Recall that  $\mathcal{P} = (2, x_1-1, \ldots, x_l-1)$  is a maximal ideal in $\mathbb{Z}_4[x_1, x_2, \ldots, x_l]$ with residue field $\mathbb{F}_2$.) By Lemma \ref{Rislocal}, the ring $R$ is a local ring, and therefore its quotient $S$ is also a local ring whose maximal ideal corresponds to $\mathcal{P}/\mathcal{J}$ with residue field  $\mathbb{F}_2$.
This gives a short exact sequence:
\[0 \rightarrow \mathcal{P}/\mathcal{J} \rightarrow S \rightarrow \mathbb{F}_2 \rightarrow 0.\]
Since $\mathbb{F}_2^\times = \{ 1\}$, and units in a local ring are the elements that are in the complement of the maximal ideal (see Lemma \ref{keylemma}),  it is clear from the above short exact sequence that any unit in $S$ is of the form $\overline{1 + \eta}$ where $\eta$ belongs to $\mathcal{P}$. With this characterization of the units at hand, we are now ready to characterize the $\Delta_2$ quotients of $R$. 

Suppose that a quotient $S$ of $R$ is a $\Delta_2$-ring.  Consider any $\eta$ in $\mathcal{P}$. Since $S$ is  $\Delta_2$-ring, the unit $\overline{1 + \eta}$ in $S$ must square to 1. Note that $(\overline{1 + \eta})^2 = \overline{1 +\eta^2 +2 \eta} =\overline{1 +\eta(\eta +2})$. This element is equal to  $\overline{1}$ in $S$ if and only if $\eta(\eta +2)$ belongs to $\mathcal{J}$, showing that $( \eta(\eta +2) \colon \eta \in \mathcal{P} ) \subseteq \mathcal{J}$.

Conversely, suppose $\mathcal{J}$ contains $( \eta(\eta +2) \colon \eta \in \mathcal{P} )$. Pick any unit in $S$. By the above discussion, the chosen unit has to be of the form $\overline{1 + \eta}$, where $\eta$ belongs to $\mathcal{P}$. We have, 
$(\overline{1 + \eta})^2 = \overline{1 +\eta^2 +2 \eta} =\overline{1 +\eta(\eta +2})$. This last expression is $\overline{1}$ in $S$ because $\eta(\eta +2)$ belongs to $\mathcal{J}$. This shows that $S$ is a $\Delta_2$-ring.
\end{proof}

In characteristic 8, we have the following result similar to the previous one with almost identical proof. We use the fact that $2^3 = 0$ in $\mathbb{Z}_8$.

\begin{prop}
Let $R = \mathbb{Z}_8[x_1, x_2, \ldots, x_l]/(x_1^2-1, x_2^2-1, \ldots, x_l^2-1)$ and let $\mathcal{P} = (2, x_1-1, x_2-1, \ldots, x_l-1)$ be the maximal ideal of $R$. The ring $S$, a quotient of $R$, is a $\Delta_2$-ring if and only if 
$S = R/J$ where $J$ contains $( \eta(\eta + 2) \colon \eta \in \mathcal{P} )$. (Note: $\mathcal{P}$  is the unique prime ideal of $R$.)
\end{prop}

Packing all the above lemmas and propositions, we get the following theorem that gives a complete characterization of finite commutative rings that are $\Delta_2$.

\begin{thm}
Let $R$ be a finite commutative ring that is $\Delta_2$. Then the characteristic of $R$ is a divisor of $24$. In each characteristic $n$ that divides $24$, the finite commutative rings that are $\Delta_2$ are all quotients (with characteristic $n$) of the rings shown in Table \ref{finitecommutativering}. 

\begin{table}[!h] 
    \centering
    \begin{tabular}{c|c}
        $\text{char}(R)$ & $R$  \\
         \hline\\
         $2$ & $\prod \; \mathbb{F}_2[x_1, x_2, \ldots, x_{l}]/(x_1^2, x_2^2, \ldots, x_{l}^2)$ \\ \\
         $3$ & $ \prod \; \mathbb{F}_3$ \\ \\
         $6$ & $\prod \; \mathbb{F}_3  \times \prod (\mathbb{F}_2[x_1, x_2, \ldots, x_{l}]/(x_1^2, x_2^2, \ldots, x_{l}^2)) $ \\ \\
         $4$ &  $ \prod \; \mathbb{Z}_4[x_1, x_2, \ldots, x_{l}]/(  r(r+2)  \colon r \in \mathcal{P} )  $ \\ \\
         $8$ &  $ \prod \; \mathbb{Z}_8[x_1, x_2, \ldots, x_{l}]/(  r(r+2)  \colon r \in \mathcal{P} )  $ \\ \\
         $12$ &  $ \prod \; \mathbb{F}_3  \times \prod \;  \mathbb{Z}_4[x_1, x_2, \ldots, x_{l}]/(  r(r+2)  \colon r \in \mathcal{P} )$ \\ \\
         $24$ &  $ \prod \; \mathbb{F}_3  \times \prod \; \mathbb{Z}_8[x_1, x_2, \ldots, x_{l}]/(  r(r+2)  \colon r \in \mathcal{P})  $ \\
    \end{tabular}\\
    \vskip 5mm
   \caption{Classification of finite commutative $\Delta_2$-rings}
   \label{finitecommutativering}
\end{table}
\end{thm}

\noindent
\textbf{Note:} In the above table, $\mathcal{P}$ is the  prime ideal $(2, x_1 -1, x_2-1, \ldots, x_{l} -1)$ in the corresponding polynomial ring. 

\begin{thm}
 
Let $R$ be a finite commutative ring, and let $p$ be an odd prime. Then, $R$ is a  $\Delta_p$-ring if and only if 
$p$ is  a Mersenne prime and  $R$ is isomorphic to  $(\mathbb{F}_{2})^a \times (\mathbb{F}_{p+1})^b$ for some $a$ and $b \ge 0$.
\end{thm} 

\begin{proof}
First, note that $\Delta_p$, for $p$ odd, can occur only when $p$ is a Mersenne prime. So, let $p=2^l-1$ for some $l$.  As before, we may assume without loss of generality that $R$ is a local Artin ring. Let $n$ be the characteristic of $R$.  Then $\mathbb{Z}_n$ is a $\Delta_p$-ring. That is $u^p=1$ for all $u$ in $\mathbb{Z}_n^\times$. If $n \ne 2$, $\mathbb{Z}_n^\times$ will have an element of order $2$, which is impossible in a $\Delta_p$-ring. This shows that $n =2$.   By Lemma \ref{keylemma}, $R$ then has to be to be a quotient of $\mathbb{F}_2[C_{p}^l]$.
By \ref{commutativeSS}, a  quotient of $\mathbb{F}_2[C_{p}^l]$ must be a product of fields of characteristic $2$. Since the multiplicative group must be a elementary abelian $p$-group, the finite fields can be either $\field_2$ or $\field_{p+1}$ with $p$ Mersenne. This shows that $R$ must have the form stated in the theorem. Conversely, it is easy to see that these rings have elementary abelian $2$-groups as their unit groups.
\end{proof}

\section{Modular group algebras: $\mathbb{Z}_nG$} \label{mgr}

\begin{thm}
    The group algebra $\zed_nG$ is $\Delta_2$ if and only if $n \in \{2,3,6\}$ and 
    $G = C_2^r$ or $n \in \{4,12\}$ and $G = C_2$.
\end{thm}

\begin{proof}
 Let $\zed_nG$ be a $\Delta_2$-ring. Then $\zed_n$ must be $\Delta_2$, as it is
 a subring of $\zed_nG$. Then by Lemma \ref{24}, we know that $n$ divides 24.
 Note that every element of $G$ is a unit of $\zed_nG$, thus we have $g^2 = e$ for any  $g \in G$.  Proposition \ref{elemabelian} implies that $G \cong C_2^r$,  an elementary abelian $2$-group.

 Let $C_2 = \{e, \sigma\}$, and denote the 
 $i$th generator of $C_2^r$ by $\sigma_i$.  We proceed by checking each divisor of 24. We use the fact that the sum of a unit and a nilpotent element in any commutative ring is a unit. 
 \begin{itemize}
     \item $\zed_2C_2^r$  and  $\zed_3C_2^r$ are $\Delta_2$-rings from \cite[Theorem 2.3]{CLY}.
     \item Note that $\zed_8C_2^r$ is not $\Delta_2$, as $u = e + 2\sigma$ is a unit,
     but $u^2 = 5e+4\sigma \ne e$.
     \item Similarly, $\zed_{24}C_2^r$ is not $\Delta_2$, as $u = e+ 6\sigma$ is a unit,
     but $u^2 = 13e + 12\sigma \ne e$. 
     \item $\zed_6C_2^r$ is isomorphic to $\zed_2C_2^r \times \zed_3C_2^r$ via the 
     Chinese Remainder Theorem, thus because each of $\zed_2C_2^r$ and $\zed_3C_2^r$
     are $\Delta_2$, we know that $\zed_6C_2^r$ is $\Delta_2$. 
     \item To show that $\zed_4C_2$ is $\Delta_2$, consider the augmentation 
     map $\epsilon: \zed_4C_2 \to \zed_4$ defined by $\sum \alpha_g g \mapsto \sum \alpha_g$. 
     Then $\epsilon$ is a ring homomorphism, so it must map unit elements of $\zed_4C_2$ to 
     unit elements of $\zed_4$. The only unit elements of $\zed_4$ are 
     $\pm 1$, and it can be verified that each of the elements in $\{e, \sigma, 3\sigma, e + 2\sigma,
     2e + \sigma, 2e + 3\sigma, 3e, 3e + 2\sigma\}$ are units of $\zed_4C_2$ satisfying
     $u^2 = e$. 
     However, 
     for $\zed_4C_2^r$ where $r > 1$, note that $\sigma_1 + \sigma_2$ is a nilpotent 
     as $(\sigma_1 + \sigma_2)^4 = 0$. Then $u = e + \sigma_1 + \sigma_2$ is a unit in 
     $\zed_4C_2^r$, but $u^2 = 3e + 2\sigma_1 + 2\sigma_2 + 2\sigma_1\sigma_2 \ne e$. 
     \item We know by the Chinese Remainder Theorem that for all $r \ge 1$,  $\zed_{12}C_2^r \cong \zed_4C_2^r \times 
     \zed_3C_2^r$ (because $4$ and $3$ are relatively prime). When $r=1$, we know from previous cases that both $\zed_4C_2$ and $\zed_3C_2$ are $\Delta_2$. It follows that  that  $\zed_{12}C_2$ is $\Delta_2$ as well.
     However, when $r > 1$,  $\zed_4C_2^r$ was shown to be not $\Delta_2$. Therefore, $\zed_{12}C_2^r$ also can't be $\Delta_2$.
 \end{itemize}
\end{proof}

\begin{thm}
 Let $p$ be an odd prime. The group ring
$\mathbb{Z}_nG$ ($G$ abelian) is a $\Delta_p$-ring if and only if $n=2$, $G$ is an elementary abelian $p$-group, and $p$ is some Mersenne prime.
\end{thm}

\begin{proof}
Let $\mathbb{Z}_nG$ with $G$ abelian be a $\Delta_p$-ring. 
Since $G$ is abelian  and  $\mathbb{Z}_nG$ is $\Delta_p$,  every non-trivial element in $G$ has order $p$.   This means $G$ has to be an elementary abelian $p$-group. Since the subring $\mathbb{Z}_n$ must also be a  $\Delta_p$-ring,  the characteristic $n$ must be $2$, because otherwise, we will have $-1 (\ne 1)$, a unit element in the ring of order $2$, impossible in a $\Delta_p$-ring with $p$ odd. So our group ring is $\field_2C_p^r$. This is $\Delta_p$ if and only if $p$ is Mersenne; \cite{CLY}. 

\end{proof}

The next proposition gives an ideal-theoretic explanation for why $\mathbb{Z}_4[C_2]$ is $\Delta_2$, but  $\mathbb{Z}_4[C_2^l]$ is not  $\Delta_2$ when $l > 1$. In the ring $\mathbb{Z}_4[x_1, \ldots, x_l]$, consider the following ideals. 

\begin{itemize}
    \item $P_l := ( 2, x_1-1,x_2-1, \ldots, x_l-1 )$
    \item $J_l:= ( n(n+2) : n \in P_l )$
    \item  $I_l := ( x_1^2-1, \ldots, x_l^2-1 )$
\end{itemize}

\begin{prop} \label{idealrelationship} Let $l$ be a positive integer. The following are equivalent.
\begin{enumerate}
    \item $l = 1$.
    \item $\mathbb{Z}_4[C_2^l]$ is a $\Delta_2$-ring.
    \item $I_l = J_l$.
\end{enumerate}
\end{prop}

\begin{proof}
 We have already seen the equivalence of (1) and (2). We will show that (1) and (3) are equivalent.
 To this end, we have to show that, in the ring $\mathbb{Z}_4[x]$, $( x^2-1 ) = ( n(n+2) \colon \colon n \in P )$ where $P = ( 2,x-1 )$. The inclusion $( x^2-1 ) \subseteq ( n(n+2) \colon  n \in P )$  is obvious because $x^2 -1 = (x-1)(x+1) = (x-1)((x-1)+2)$. For the other inclusion, it is enough to show that for all $n$ in $P$, the element $n(n+2)$ is $0$ in the quotient ring $\mathbb{Z}_4[x]/(x^2-1)$. Note that all multiples of 4 will be 0 and $x^2 =1$ in the quotient ring. Keeping this in mind, consider an arbitrary element $n(n+2)$ in this quotient ring, where $n = 2(c+bx) + (x-1)(a+bx)$ is in $P$. Then we have the following equations.
 \begin{eqnarray*}
 n(n+ 2) & =&  n^2 +2n \\
 & = & (2(c+dx) + (x-1)(a+bx))^2 + 2(2(c+dx) + (x-1)(a+bx))\\
   & = & (x-1)^2 (a+bx)^2 + 2(x-1)(a+bx) \\
   & = & 2(1-x)(a^2+b^2) + 2(x-1)(a+bx) \\
   & = & 2(1-x)(a^2+b^2-a -bx)
 \end{eqnarray*}
The last expression is directly seen to be zero in our quotient ring for any choice of $a$ and $b$ in $\mathbb{Z}_4$. For $l > 1$, we claim that $I_l \subsetneq  J_l$. It is clear that $I_l \subseteq  J_l$. To see that the inclusion is strict, note that, for $l > 1$,  from the above results we have $\mathbb{Z}_4[x_1, \ldots, x_l]/J_l$ is a $\Delta_2$-ring but  $\mathbb{Z}_4[x_1, \ldots, x_l]/I_l$ is not a $\Delta_2$-ring. This shows that $I_l \subsetneq  J_l$ for $l > 1$. 
\end{proof}

\section{Path algebras: $kQ$} \label{pa}
    
    Let $Q$ be a quiver (a directed, not necessarily simple graph), and $k$ be a field. 
    We define the path algebra $kQ$ to be a vector space 
    over $k$ with basis given by paths in $Q$, including the trivial paths of length $0$ starting 
    and ending at the same vertex (we will denote the trivial path of length $0$ starting and ending
    at a vertex $i$ by $e_i$). For any two paths $p,q$, we define the multiplication $pq$ to be 
    the concatenation of $p$ and $q$ if $t(q) = s(p)$, and 0 otherwise; here $t(q)$ is the tail of $q$ and $s(p)$ is the head of $p$. 
    
    The existence of an identity element is guaranteed in any finite quiver. Moreover, the path algebra $kQ$ of a quiver $Q$ will be finite-dimensional if and only if $Q$ is acyclic. So, we will assume that our quivers are finite and acyclic.
    
\begin{lem}[Karthika-Viji \cite{KV}] \label{identitykq}
   Let $k$ be a field, and $Q$ be a finite acyclic quiver. Then the identity element of 
    the path algebra $kQ$ is given by $\sum_{i \in V(Q)} e_i$. We denote this element by $e$.
\end{lem}


   The following theorem gives a useful characterization of units in path algebras.

\begin{thm}[Karthika-Viji \cite{KV}]\label{unitkq}
    Let $k$ be a field, and $Q$ be a finite acyclic quiver. 
    Then an element $a \in kQ$ is a unit if and only if the coefficient of $e_i$ is nonzero for all vertices $i \in Q$.
\end{thm}

The following two results complete the classification of all $\Delta_p$-path algebras. 

\begin{thm}
    Let $k$ be a field, and $Q$ be a finite acyclic quiver. Then $kQ$ is $\Delta_2$ if and only if
    $k = \field_3$ and $Q$ contains no edges, or $k = \field_2$ and $Q$ has no directed paths of length $2$.
\end{thm}

\begin{proof}
    Let $k = \field_3$ and $Q$ be a quiver containing no nontrivial path. Then by Theorem \ref{unitkq}, 
    the only unit in $kQ$ is the identity $e = \sum e_i$. Then we see that $e^2 = e$, as 
    $e$ is the identity. 
    
    Let $k = \field_3$ and $Q$ be a quiver containing a directed edge $p$. By Theorem \ref{unitkq}, $e + p$ is a unit, but $(e + p)^2 = e^2 + ep + pe + p^2 = e + p + p + 0 = e + 2p \ne e$. Thus $kQ$ is not $\Delta_p$.
    
    Let $k = \field_2$ and $Q$ be a quiver containing a directed path $\beta\alpha$ comprised of edges $\alpha, \beta$. Then 
    by Theorem \ref{unitkq}, $e + \alpha + \beta$ is a unit, but $(e+\alpha+\beta)^2 = e^2 + 2\alpha + 2\beta + \beta\alpha \ne e$.
    
    If $k = \field_2$ and $Q$ is a quiver with no directed paths of length two, then for any two paths $\alpha, \beta$,
    we have $\alpha\beta = 0$. Let $p$ be the sum of any paths. Then $e + p$ is a unit, but
    $(e+p)^2 = e + 2p = e$ as desired. 
\end{proof}
\begin{example}  Here is an example of a quiver on $5$ vertices for which the path algebra over $\field_2$ or $\field_3$ is a $\Delta_2$-ring.
\[
\begin{tikzcd}
\bullet &&& \bullet &&& \bullet &&& \bullet && \bullet
\arrow[from=1-7, to=1-10]
\arrow[from=1-12, to=1-10]
\arrow[from=1-7, to=1-4]
\arrow[from=1-1, to=1-4]
\end{tikzcd}
\]
\end{example}
Note that this quiver has no directed path of length $2$ or more.
\begin{thm}
    Let $p$ be an odd prime, $k$ be a field, and $Q$ be a finite acyclic quiver. The path algebra 
    $kQ$ is $\Delta_p$ if and only if $Q$ is trivial and $k = \field_2$ or $k = 
    \field_{p+1}$ where $p$ is a Mersenne prime. 
\end{thm}

\begin{proof}
Let $kQ$ be a $\Delta_p$-ring for some odd prime $p$.  Note that $kQ$ is a $k$-algebra. Then by Corollary \ref{deltapfield}, $k$ must be $\field_2$ or $k = \field_{p+1}$   where $p$ is Mersenne.   We claim that $Q$ cannot have any directed edges. Suppose to the contrary, $Q$ has some directed edge $\alpha$. Consider the element $u = e + \alpha$. Then 
 $u$ is a unit by Lemma $\ref{unitkq}$, and $u^2 = (e + \alpha)^2 = e + 2\alpha$. 
 Since the characteristic of $k$ is 2, $2\alpha = 0$, thus $u^2 = e$ and we see that  $u$ has order 2. This contradicts the fact that $kQ^\times = C_p^r$, as every element of $kQ^\times$ should have order $p$.  This proves one direction. The other direction is obvious because the path algebra $kQ$ of a trivial quiver is a product of copies of $k$, and the direct product of fields stated in the theorem is a $\Delta_p$-ring.
\end{proof}

\begingroup
\raggedright

\bibliographystyle{plain}
\bibliography{citations}
\endgroup

\end{document}